\documentclass[11pt]{amsart}
\usepackage{latexsym,amsmath,amssymb}

\title{\protect{Sobolev spaces, Lebesgue points and maximal functions}}

\author{Piotr Haj\l{}asz, Zhuomin Liu}

\address{Department of Mathematics, University of Pittsburgh, 301
  Thackeray Hall, Pittsburgh, PA 15260, USA, {\tt hajlasz@pitt.edu}}

\address{Department of Mathematics, University of Pittsburgh, 301
  Thackeray Hall, Pittsburgh, PA 15260, USA, {\tt liuzhuomin@hotmail.com}}

\thanks{P.H. was supported by NSF grant DMS-1161425.}

plus 6pt minus 12pt
\def\eps{\varepsilon}

\def\vi{\varphi}

\def\M{{\mathcal M}}
\def\H{{\mathcal H}}

\def\S{{\mathcal S}}
\def\I{{\mathcal I}}
\def\T{{\mathcal T}}

\def\rss{{\vert_{_{_{\!\!-\!\!-\!}}}}}
\newtheorem{theorem}{Theorem}
\newtheorem{lemma}[theorem]{Lemma}

\newtheorem{proposition}[theorem]{Proposition}

\theoremstyle{definition}

\newtheorem{example}[theorem]{Example}
\newtheorem{question}[theorem]{Question}

\newcommand{\barint}{
\rule[.036in]{.12in}{.009in}\kern-.16in \displaystyle\int }

\newcommand{\barcal}{\mbox{$ \rule[.036in]{.11in}{.007in}\kern-.128in\int $}}

\newcommand{\bbbr}{\mathbb R}

\def\bbbc{{\mathchoice {\setbox0=\hbox{$\displaystyle\rm C$}\hbox{\hbox
to0pt{\kern0.4\wd0\vrule height0.9\ht0\hss}\box0}}
{\setbox0=\hbox{$\textstyle\rm C$}\hbox{\hbox
to0pt{\kern0.4\wd0\vrule height0.9\ht0\hss}\box0}}
{\setbox0=\hbox{$\scriptstyle\rm C$}\hbox{\hbox
to0pt{\kern0.4\wd0\vrule height0.9\ht0\hss}\box0}}
{\setbox0=\hbox{$\scriptscriptstyle\rm C$}\hbox{\hbox
to0pt{\kern0.4\wd0\vrule height0.9\ht0\hss}\box0}}}}

\def\bbbq{{\mathchoice {\setbox0=\hbox{$\displaystyle\rm Q$}\hbox{\raise
0.15\ht0\hbox to0pt{\kern0.4\wd0\vrule height0.8\ht0\hss}\box0}}
{\setbox0=\hbox{$\textstyle\rm Q$}\hbox{\raise 0.15\ht0\hbox
to0pt{\kern0.4\wd0\vrule height0.8\ht0\hss}\box0}}
{\setbox0=\hbox{$\scriptstyle\rm Q$}\hbox{\raise 0.15\ht0\hbox
to0pt{\kern0.4\wd0\vrule height0.7\ht0\hss}\box0}}
{\setbox0=\hbox{$\scriptscriptstyle\rm Q$}\hbox{\raise
0.15\ht0\hbox to0pt{\kern0.4\wd0\vrule height0.7\ht0\hss}\box0}}}}

\def\bbbz{{\mathchoice {\hbox{$\sf\textstyle Z\kern-0.4em Z$}}
{\hbox{$\sf\textstyle Z\kern-0.4em Z$}} {\hbox{$\sf\scriptstyle
Z\kern-0.3em Z$}} {\hbox{$\sf\scriptscriptstyle Z\kern-0.2em
Z$}}}}

\def\mvint_#1{\mathchoice
          {\mathop{\vrule width 6pt height 3 pt depth -2.5pt
                  \kern -8pt \intop}\nolimits_{\kern -3pt #1}}%
          {\mathop{\vrule width 5pt height 3 pt depth -2.6pt
                  \kern -6pt \intop}\nolimits_{#1}}%
          {\mathop{\vrule width 5pt height 3 pt depth -2.6pt
                  \kern -6pt \intop}\nolimits_{#1}}%
          {\mathop{\vrule width 5pt height 3 pt depth -2.6pt
                  \kern -6pt \intop}\nolimits_{#1}}}

\numberwithin{theorem}{section} \numberwithin{equation}{section}

\begin{document}




\begin{abstract}
In this note we study boundedness of a large class of maximal operators 
in Sobolev spaces that includes the
spherical maximal operator. We also study the size of the set of Lebesgue 
points with respect to convergence associated with such maximal operators.  
\end{abstract}

\maketitle

\begin{center}
\it To Professor Bogdan Bojarski
\end{center}

\section{Introduction}

Kinnunen \cite{kinnunen} proved that the Hardy-Littlewood maximal operator is bounded in the Sobolev space $W^{1,p}(\bbbr^n)$ for $1<p\leq\infty$. Actually a slightly more general result is true \cite[Theorem~1]{hajlaszo}, \cite[Theorem~2.1]{kinnunens}. We say that an operator $A$ defined on a linear space of measurable functions on $\bbbr^n$ is {\em sub-linear} if $Au\geq 0$ a.e. and $A(u+v)\leq Au+Av$ a.e. We say that it {\em commutes with translations} if $A(u_y)(x)=(Au)_y(x)$, where $u_y(x)=u(x-y)$.
\begin{proposition}
If $A:L^p(\bbbr^n)\to L^p(\bbbr^n)$, $1<p<\infty$ is bounded, sub-linear and commutes with translations, then
$A:W^{1,p}(\bbbr^n)\to W^{1,p}(\bbbr^n)$ is bounded.
\end{proposition}
The proof easily follows from a characterization of the Sobolev space in terms of difference quotients, \cite[Section~7.11]{gilbarg}. 
Indeed, if $u\in W^{1,p}(\bbbr^n)$, $1<p<\infty$, then
$
\Vert u_y-u\Vert_p\leq \Vert\nabla u\Vert_p|y|
$
and if $u\in L^p$, $1<p<\infty$ satisfies
$
\Vert u_y-u\Vert_p\leq C|y|
$
for all $y\in\bbbr^n$,
then $u\in W^{1,p}$ and $\Vert \nabla u\Vert_p\leq C$.
Hence
the result follows from the following estimates
\begin{eqnarray*}
\Vert (Au)_y-Au\Vert_p
& = &
\Vert A(u_y)-Au\Vert_p\leq
\Vert A(u-u_y)\Vert_{p} +
\Vert A(u_y-u)\Vert_p\leq \\
& \leq &
C\Vert u_y-u\Vert_p \leq C\Vert \nabla u\Vert_p |y|.
\end{eqnarray*}
$\phantom{a}$
\hfill $\Box$

According to the theorem of Bourgain and Stein \cite{bourgain}, \cite{stein}, (see also \cite{grafakos}, \cite{steinbook}) the spherical maximal operator
$$
\S u(x) = \sup_{t>0}\, \barint_{S^{n-1}(x,t)} u(z)\, d\H^{n-1}(z)
$$
is bounded in $L^p(\bbbr^n)$ when $n\geq 2$ and $p>n/(n-1)$. 
Hence $\S:W^{1,p}(\bbbr^n)\to W^{1,p}(\bbbr^n)$ is bounded for $n\geq 2$ and $p>n/(n-1)$. 
Since for $1<p<n$ we have $W^{1,p}(\bbbr^n)\subset L^{p^*}(\bbbr^n)$, where 
$p^*=np/(n-p)>n/(n-1)$ we can go with the exponent $p$ below $n/(n-1)$ and conclude that 
$\S:W^{1,p}\to L^{p^*}$ is bounded for all $1<p\leq n/(n-1)$. 

One can easily construct a function $u\in L^p$, $1\leq p\leq n/(n-1)$ such that $\S u\equiv\infty$, see
\cite{grafakos}, \cite{steinbook}.
Hence $\S$ is not bounded in $L^p$ for $1\leq p\leq n/(n-1)$. 
In this note among other facts we will show that for $u\in W^{1,p}$, $1<p<n$,
the function $\S u$ has growth properties of a $\dot{W}^{1,p}$ function, see \eqref{p-pstar} and Theorem~\ref{spherical-generalized}, 
and hence
examples for the lack of boundedness of $\S$ in $L^p$, $1<p\leq n/(n-1)$ {\em cannot} be easily generalized to the Sobolev case
(for the case $p=1$ see, however Example~\ref{E1} below). 
Here $\dot{W}^{1,p}$, $1\leq p<n$ stands for the homogeneous Sobolev space which is defined by
$$
\dot{W}^{1,p}(\bbbr^n)=\{u\in L^{p^*}:\, \nabla u\in L^p\}.
$$
Clearly $W^{1,p}\subset \dot{W}^{1,p}\subset L^{p^*}$.

More precisely, we will provide an elementary proof of the fact that for a wide range of maximal operators $\S_\mu$,
that includes the spherical maximal operator, there is another operator $\T$ such that $\S_\mu\leq C\T$ and 
$\T:W^{1,p}\to \dot{W}^{1,p}$ is bounded for all $1<p<n$.
This gives some evidence that the answer to \cite[Question~2]{hajlaszo} might be in the positive, see also Question~\ref{Q1} below.

In particular $\S_\mu:W^{1,p}\to L^{p^*}$ is bounded for all $1<p<n$. 
The class of maximal operators $\S_\mu$ includes the spherical maximal operator,
but also the maximal operator where we take averages over boundaries of cubes. Note that such a maximal operator is 
{\em not} bounded in $L^p$ for any $p$ and hence our argument cannot involve
the Bourgain-Stein result, so the argument has to be substantially different from that used in \cite{hajlaszo}.
It also shows that the boundedness properties
of the spherical maximal operator in the Sobolev setting are not necessarily based on
the boundedness of the spherical maximal operator in $L^p$, another evidence that we could go with the exponent below $n/(n-1)$.

As an application of our result we will also provide a generalization of the Federer-Ziemer theorem \cite{federerz} about the size of the set of Lebesgue points of a Sobolev function, Theorem~\ref{FZ-leb}.

Let $\mu$ be a probability Borel measure on $\bbbr^n$. We define the rescaled measure $\mu_t$, $t>0$ by
$$
\mu_t(E)=\mu(E/t),
\qquad
E/t=\{x/t:\, x\in E\}
$$
Clearly
$$
\int_{\bbbr^n} u(z)\, d\mu_t(z) =
\int_{\bbbr^n} u(tz)\, d\mu(z).
$$
If $\sigma=(n\omega_n)^{-1}\H^{n-1}\rss S^{n-1}(0,1)$ is the normalized Hausdorff measure on the unit sphere, then
$\sigma_t=(n\omega_n t^{n-1})^{-1}\H^{n-1}\rss S^{n-1}(0,t)$ is the normalized Hausdorff measure on 
the sphere of radius $t$, and if $\nu$ is the normalized Lebesgue measure on the unit ball, then $\nu_t$ is
the normalized Lebesgue measure on the ball of radius $t$.
Here $\omega_n$ stands for the volume of the unit ball in $\bbbr^n$ and hence $n\omega_n$ is the volume of the unit sphere.

With the measure $\mu$ we associate the following maximal operator
$$
\S_\mu u(x) = \sup_{t>0}\int_{\bbbr^n} |u(z+x)|\, d\mu_t(z) 
$$
If $\mu=\sigma$, we obtain the spherical maximal operator 
$\S_\sigma u=\S u$ and if
$\nu$ is the normalized Lebesgue measure on the unit ball we obtain the Hardy-Littlewood maximal operator
$$
\S_\nu u(x) = \M u(x) = \sup_{t>0}\barint_{B(x,t)} |u(z)|\, dz.
$$
We say that a probability measure $\mu$ 
is {\em spherical-like} if it is supported on a bounded set, 
\begin{equation}
\label{sl}
{\rm supp}\,\mu\subset\overline{B}(0,R)
\end{equation}
and satisfies the estimate
\begin{equation}
\label{spherical-like}
\sup_{x\in \bbbr^n, r>0} \frac{\mu(B(x,r))}{r^{n-1}}=M<\infty.
\end{equation}
Note that the measures $\sigma$ and $\nu$ have this property.
Moreover the normalized Hausdorff measure on the boundary of a cube or on any other 
compact $(n-1)$-dimensional Ahlfors regular set has this property.

Importance of this condition stems from the following
beautiful result due to Meyers and Ziemer, \cite{meyersz}, \cite[Lemma~4.9.1]{ziemer}.
\begin{lemma}
\label{meyers-ziemer}
If $\mu$ is a positive Radon measure on $\bbbr^n$,
then there is a constant $C>0$ such that
\begin{equation}
\label{poincare}
\int_{\bbbr^n} |u|\, d\mu \leq C\int_{\bbbr^n} |\nabla u|\, dx
\quad
\mbox{for all $u\in C_0^\infty(\bbbr^n)$}
\end{equation}
if and only if the condition \eqref{spherical-like} is satisfied. Moreover \eqref{poincare} holds with $C=c(n)M$.
\end{lemma}
In this result we do not assume that $\mu$ is a probability measure neither that it is supported on a bounded set.

The {\em Riesz potential} is defined by
$$
\I g(x)=\int_{\bbbr^n}\frac{g(z)}{|x-z|^{n-1}}\, dz.
$$
The classical Fractional Integration Theorem \cite[Theorem~2.8.4]{ziemer}
asserts that
$$
\I: L^p\to L^{p^*},
\quad 
1<p<n
$$
is a bounded operator. Actually, this is a consequence of a much deeper result 
\cite[Chapter~5]{stein2} (see also \cite{hajlaszl}) which states that
\begin{equation}
\label{riesz-sob}
\I :L^p\to \dot{W}^{1,p},
\quad
1<p<n
\end{equation}
is a bounded operator.
Let 
$$
\T u=\M u+\I |\nabla u|.
$$ 
According to \eqref{riesz-sob} and Kinnunen's theorem about boundedness of the maximal operator we have that
\begin{equation}
\label{p-pstar}
\T :W^{1,p}\to \dot{W}^{1,p},
\quad
1<p<n
\end{equation}
is bounded. In particular
$$
\T :W^{1,p}\to L^{p^*},
\quad
1<p<n
$$
is also bounded.

Our first result reads as follows.
\begin{theorem}
\label{spherical-generalized}
Let $\mu$ be a spherical-like measure
satisfying \eqref{sl} and \eqref{spherical-like}.
Then there is a constant $C=C(n)R^{n-1}M>0$ such that
$$
\S_\mu u\leq C\T u
\quad
\mbox{everywhere}
$$
for all $u\in W^{1,p}$, $1<p<n$. In particular
\begin{equation}
\label{without-Bourgain}
\S_\mu:W^{1,p}\to L^{p^*},
\quad
1<p<n
\end{equation} 
is bounded.
\end{theorem}
As we mentioned earlier, the result applies to the spherical maximal operator, but 
also to the maximal operator where we take averages over boundaries of cubes.
The conclusion \eqref{without-Bourgain} was proved in \cite{hajlaszo} 
for the spherical maximal operator with the use 
the Bourgain-Stein  theorem. However, our proof of 
is elementary and it does not involve the Bourgain-Stein result.

The operator $\S_\mu$ is bounded by the operator $\T$ which preserves Sobolev classes in the sense of \eqref{p-pstar}, but it does not necessarily imply that the operator $\S_\mu$ has a similar property. 
While the function $\S_\mu u$ has growth properties of a Sobolev function it may happen that it has high oscillations which
possibly could exclude it from being in the Sobolev space.
However, Theorem~\ref{spherical-generalized} suggests the following question.

\begin{question}
\label{Q1}
Is the spherical maximal operator bounded in the Sobolev space
$\S:W^{1,p}(\bbbr^n)\to \dot{W}^{1,p}(\bbbr^n)$ for all $n\geq 2$ and $1<p<n$?
\end{question}

Example~\ref{E1} below shows that in the case $p=1$ the corresponding question
has the negative answer.

We could formulate the question in the case of more general maximal operators, but 
one should treat first the case of the spherical maximal operator especially
that in that case many nice integral formulas and connections to singular integrals are available, see \cite{hajlaszl}.
This is why we did not state the problem in a more general form.

\begin{example}
\label{E1}
We will show that there is $u\in W^{1,1}(\bbbr^n)$ vanishing outside a compact set such that $\S u$ is not even in 
$W^{1,1}_{\rm loc}(\bbbr^n)$. 
Let $u$ be a smooth extension of
$$
|x|^{1-n}\log^{-1-(n-1)/n}(e/|x|)
$$ 
from the unit ball
to a compactly supported function. Then
$u\in W^{1,1}(\bbbr^{n})$. 
Suppose that $\S u\in W^{1,1}_{\rm loc}$. Then 
$\S u\in L^{n/(n-1)}_{\rm loc}$. On the other hand
a simple computation shows that for $|x|\leq 1$
$$
\S u(x) \geq v(x):=\mvint_{S^{n-1}(x,|x|)}|u|\, d\sigma
\geq \frac{C}{|x|^{n-1}(\log(e/|x|))^{(n-1)/n}}
$$
and clearly the right hand side is not in $L^{n/(n-1)}$ in any neighborhood
of the origin.
\hfill $\Box$
\end{example}

Federer and Ziemer \cite{federerz} proved that the set of non-Lebesgue points of 
a $p$-quasicontinuous representative of a Sobolev function $f\in W^{1,p}$ has $p$-capacity zero.

Recall that the $p$-capacity, $1<p<n$ of a set $A\subset\bbbr^n$ is defined as
$$
C_p(A)=\inf\left\{\int_{\bbbr^n}|\nabla u|^p\, dx\right\}
$$
where the infinum is taken over all
$u\in\dot{W}^{1,p}$ such that $u\geq 0$ and  $u\geq 1$ 
in an open set that contains $A$.
For basic properties of the capacity, see \cite[Section~4.7]{EG}.

A function $u$ is said to be {\em $p$-quasicontinuous} if for every $\eps>0$ there
is an open set $U\subset\bbbr^n$ such that
$$
C_p(U)<\eps
\quad
\mbox{and}
\quad
\mbox{$u|_{\bbbr^n\setminus U}$ is continuous.}
$$
Every function $u\in W^{1,p}(\bbbr^n)$, $1<p<n$ has a $p$-quasicontinuous representative, see 
\cite[Section~4.8]{EG} and any two $p$-quasicontinuous representatives of
$u\in W^{1,p}$ are equal away from a set of $p$-capacity zero \cite{kilpelainen}.
It is well known that for $u\in W^{1,p}$ the following representative defined at every point of $\bbbr^n$ is $p$-quasicontinuous
\begin{equation}
\label{precise}
u(x):=\limsup_{r\to 0}\barint_{B(x,r)} u(z)\, dz,
\end{equation}
see, \cite[Section~4.8]{EG}.
Clearly it suffices to prove the Federer-Ziemer theorem for the representative given by
\eqref{precise}.
Theorem~\ref{spherical-generalized}, or rather its proof, leads to the following generalization of the Federer-Ziemer theorem.

\begin{theorem}
\label{FZ-leb}
Let $\mu$ be a spherical-like measure. Let $u\in W^{1,p}(\bbbr^n)$,
$1<p<n$ be a $p$-quasicontinuous representative. Then there is a set $E\subset\bbbr^n$ of
$p$-capacity zero $C_p(E)=0$ such that
$$
\lim_{t\to 0}\int_{\bbbr^n} |u(z+x)-u(x)|\, d\mu_t(z)=0
\quad
\mbox{for all $x\in\bbbr^n\setminus E$.}
$$
In particular
$$
\lim_{t\to 0}\int_{\bbbr^n} u(z+x)\, d\mu_t(z)=u(x)
\quad
\mbox{for all $x\in\bbbr^n\setminus E$.}
$$
\end{theorem}
In the case in which we take averages over balls, the result is due to Federer and Ziemer \cite{federerz}, but it also covers the case of taking averages over spheres and
with respect to much more general measures.

The paper is organized as follows. In Section~\ref{2} we prove Theorem~\ref{spherical-generalized}.
In Section~\ref{3} we prove Theorem~\ref{FZ-leb} and in the final Section~\ref{4} we provide a new elementary proof of
the Meyers-Ziemer theorem, Lemma~\ref{meyers-ziemer}. The original proof was based on the boxing inequality and the co-area
formula. Following ideas from \cite{malysz} we managed to avoid the co-area formula. This might have applications
to analysis on metric spaces where related estimates have been obtained with the aid of a rather involved co-area formula,
\cite{kinnunenkst}.

Throughout the paper we adopt a convention that $C$ denotes a generic constant whose value may change in a single string of estimates.

\section{Proof of Theorem~\ref{spherical-generalized}}
\label{2}

Let $\mu$ be a spherical-like measure satisfying \eqref{sl} and \eqref{spherical-like}. 
Let $\vi\in C_0^\infty(B(0,2R))$, $\vi|_{B(0,R)}\equiv 1$,
$|\nabla\vi|\leq 2R^{-1}$ be a standard cut-off function. 
First we will prove the inequality for $x=0$. Lemma~\ref{meyers-ziemer} yields
\begin{eqnarray*}
\lefteqn{\int_{\bbbr^n} |u(z)|\, d\mu_t(z)
= 
\int_{\bbbr^n} |u(tz)|\, d\mu(z) =
\int_{\bbbr^n} |\vi(z)u(tz)|\, d\mu(z)} \\
& \leq &
c(n)M\left( \int_{\bbbr^n} |\nabla\vi(z)|\, |u(tz)|\, dz +
t\int_{\bbbr^n} |\vi(z)|\, |\nabla u(tz)|\, dz\right) \\
& \leq &
c(n)M\left( R^{-1}\int_{B(0,2R)} |u(tz)|\, dz +
t\int_{B(0,2R)} |\nabla u(tz)|\, dz\right)\\
& \leq &
c(n)MR^{n-1}\left(\barint_{B(0,2tR)} |u(z)|\, dz +
\int_{B(0,2tR)}\frac{|\nabla u(z)|}{(2tR)^{n-1}}\, dz\right)\\
& \leq &
c(n)MR^{n-1}\left(\M u(0)+\I |\nabla u|(0)\right).
\end{eqnarray*}
Thus
$$
\S_\mu u(0)=\sup_{t>0} \int_{\bbbr^n} |u(z)|\, d\mu_t(z) 
\leq c(n)MR^{n-1}\left(\M u(0)+\I |\nabla u|(0)\right).
$$
Applying the inequality to $z\mapsto u(z+x)$ we get
$$
\S_\mu u(x)\leq c(n)MR^{n-1}\left( \M u(x)+\I |\nabla u|(x)\right).
$$
The proof is complete.
\hfill $\Box$

\section{Proof of Theorem~\ref{FZ-leb}}
\label{3}

We will need the following two results, see
\cite[Sections~2.4.3 and~4.7]{EG}.
\begin{lemma}
\label{hausdorff-estimate}
Let $g\in L^1_{\rm loc}(\bbbr^n)$ and $0\leq s<n$. Then
$$
\H^s\left(\left\{ x\in\bbbr^n:\, 
\limsup_{t\to 0}\frac{1}{t^s}\int_{B(x,t)}|g(z)|\, dz>0\right\}\right)=0.
$$
\end{lemma}

\begin{lemma}
\label{H-C}
For $1<p<n$ and $E\subset\bbbr^n$ we have
$$
C_p(E)\leq C\H^{n-p}(E).
$$
\end{lemma}
In particular if $g\in L^p$, $1<p<n$, then
$$
\frac{1}{t^{n-1}}\int_{B(x,t)}|g(z)|\, dz \leq
C\left(\frac{1}{t^{n-p}}\int_{B(x,t)}|g(z)|^p\, dz\right)^{1/p},
$$
and hence
\begin{equation}
\label{bla}
C_p\left(\left\{ x\in\bbbr^n:\, 
\limsup_{t\to 0} \frac{1}{t^{n-1}} \int_{B(x,t)}|g(z)|\, dz>0\right\}\right)=0.
\end{equation}
Let $u\in W^{1,p}$ be a $p$-quasicontinuous representative. 
It follows from the proof of Theorem~\ref{spherical-generalized} that
\begin{eqnarray*}
\lefteqn{\int_{\bbbr^n} |u(z+x)-u(x)|\, d\mu_t(z)} \\
& \leq &
C\left(\barint_{B(x,2tR)}|u(z)-u(x)|\, dz +
\frac{1}{(2tR)^{n-1}}\int_{B(x,2tR)}|\nabla u(z)|\, dz\right)
\end{eqnarray*}
and it suffices to observe that the first integral on the right hand side converges to
zero outside a set of $p$-capacity zero by the Federer-Ziemer theorem, while
the second integral also converges to zero outside a set of $p$-capacity zero by
\eqref{bla}.
The last conclusion of the theorem follows from the fact that $\mu_t$ is a probability measure.
The proof is complete.
\hfill $\Box$

\section{Proof of Lemma~\ref{meyers-ziemer}}
\label{4}

The necessity of the condition \eqref{spherical-like} easily follows from
\eqref{poincare} applied to suitable cut-off functions. Hence it remains to
prove that if the condition \eqref{spherical-like} is satisfied, then 
\eqref{poincare} holds with $C=c(n)M$. First we will establish a slightly weaker inequality.
\begin{lemma}
\label{level}
Let $E\subset\bbbr^n$ be a compact set and let $0\leq v\leq 1$ be a compactly supported Lipschitz function
such that $v\equiv 1$ on $E$. Then
$$
\mu(E) \leq c(n)M\int_{\bbbr^n} |\nabla v|\, dx.
$$
\end{lemma}
{\em Proof.} Consider the function
$$
\phi(t) = \int_{\{0\leq v\leq t\}} |\nabla v|\, dx.
$$
Since the function $\phi$ is increasing, it is differentiable a.e. and
$$
\int_0^1 \phi'(t)\, dx \leq \phi(1)-\phi(0)=\phi(1) = \int_{\bbbr^n} |\nabla v|\, dx.
$$
In particular there is $s\in (0,1)$ such that
$$
\phi'(s)<2\int_{\bbbr^n} |\nabla v|\, dx.
$$
This, in turn, implies that for some $\delta>0$,
$$
\frac{\phi(s)-\phi(t)}{s-t} <2 \int_{\bbbr^n} |\nabla v|\, dx
\quad
\mbox{for all $s-\delta<t<s$.}
$$
The set 
$$
E_s=\{ x:\, v(x)\geq s\}
$$
is compact and $E$ is contained in its interior. Hence from a simple
continuity of the volume argument it follows that for each $x\in E$
there is $r_x>0$ such that
\begin{equation}
\label{mz1}
|E_s\cap B(x,r_x)|=\frac{1}{2}|B(x,r_x)|.
\end{equation}
Thus also
\begin{equation}
\label{mz1.5}
\left|\{x:\, v(x)<s\}\cap B(x,r_{x})\right| = \frac{1}{2}|B(x,r_{x})|.
\end{equation}
The balls $\{B(x,r_x)\}_{x\in E}$ form an open covering of $E$ form which we
can select a finite sub-covering. Now the Vitali covering lemma
\cite[Section~1.5]{EG} implies that we can find a finite number of 
pairwise disjoint balls 
$\{B(x_i,r_{x_i})\}_{i=1}^N$ such that 
$$
E\subset \bigcup_{i=1}^N B(x_i,5r_{x_i}).
$$
Hence
\begin{equation}
\label{mz2}
\mu(E)
 \leq 
\sum_{i=1}^N \mu(B(x_i,5r_{x_i})) 
\leq 5^{n-1}M\sum_{i=1}^N r_{x_i}^{n-1}\, .
\end{equation}
Since the Lebesgue measure of the set
\begin{equation}
\label{mz3}
\{x:\, t<v(x)<s\}
\end{equation}
converges to zero as $t\to s^-$, there is $t<s$ such that the measure of the set
\eqref{mz3} is less than
$$
\min_{i\in \{1,2,\ldots,N\}} \frac{1}{4} |B(x_i,r_{x_i})|.
$$
We can also assume that $t>s-\delta$. Hence 
\eqref{mz1} and \eqref{mz1.5} yield
$$
|\{x:\, v(x)\geq s\}\cap B(x_i,r_{x_i})| = \frac{1}{2}|B(x_i,r_{x_i})|,
$$
$$
|\{x:\, v(x)\leq t\}\cap B(x_i,r_{x_i})|\geq \frac{1}{4} |B(x_i,r_{x_i})|
$$
for all $i=1,2,\ldots,N$. Now consider the truncation of $v$ between the levels 
$t$ and $s$, i.e. consider the function
$$
w(x)=
\left\{
\begin{array}{ccc}
s-t    & {\rm if} & v(x)\geq s,\\
v(x)-t & {\rm if} & t\leq v(x)\leq s,\\
0      & {\rm if} & v(x)\leq t.
\end{array}
\right.
$$
Observe that $w=s-t$ on a subset of $B(x_i,r_{x_i})$
of measure $\frac{1}{2}|B(x_i,r_{x_i})|$
and $w=0$ on a subset of $B(x_i,r_{x_i})$ of measure at least
$\frac{1}{4}|B(x_i,r_{x_i})|$.
Hence for any $c\in\bbbr$
$$
|w-c|\geq \frac{s-t}{2}
$$
on a subset of $B(x_i,r_{x_i})$ of measure at least
$\frac{1}{4}|B(x_i,r_{x_i})|$. In particular the
Poincar\'e inequality yields
$$
\frac{s-t}{8} \leq \barint_{B(x_i,r_{x_i})} |w-w_{B(x_i,r_{x_i})}|\, dx
\leq c(n)r_{x_i}\barint_{B(x_i,r_{x_i})}|\nabla w|\, dx.
$$
Here the barred integral means the integral divided by the volume of the ball
and $w_B$ is the integral average of $w$. Hence
$$
r_{x_i}^{n-1} \leq \frac{c(n)}{s-t}\int_{B(x_i,r_{x_i})} |\nabla w|\, dx
= \frac{c(n)}{s-t} \int_{B(x_i,r_{x_i})\cap \{t<v\leq s\}}|\nabla v|\, dx.
$$
Thus \eqref{mz2} and the fact that the balls $B(x_i,r_{x_i})$ are pairwise disjoint yield
\begin{eqnarray*}
\mu(E)
& \leq &
\frac{c(n)M}{s-t} \sum_{i=1}^N \int_{B(x_i,r_{x_i})\cap \{t<v\leq s\}}|\nabla v|\, dx\\
& \leq & 
\frac{c(n)M}{s-t} \int_{\{t<v\leq s\}} |\nabla v|\, dx \\
& = &
c(n)M \frac{\phi(s)-\phi(t)}{s-t} \leq c(n)M \int_{\bbbr^n} |\nabla v|\, dx.
\end{eqnarray*}
This completes the proof of Lemma~\ref{level}.
\hfill $\Box$

Now we can complete the proof of Lemma~\ref{meyers-ziemer} using the celebrated Maz'ya truncation argument
\cite{mazya} (see also \cite{hajlasz} for an expository article on that topic).
Let 
$$
u_k=
\left\{
\begin{array}{ccc}
2^{k-2}    & {\rm if} & |u|\geq 2^{k-1},\\
|u|-2^{k-2} & {\rm if} & 2^{k-2}\leq |u|\leq 2^{k-1},\\
0      & {\rm if} & |u|\leq 2^{k-2}.
\end{array}
\right.
$$
In other words $u_k$ is the truncation of $|u|$ between the levels
$2^{k-2}$ and $2^{k-1}$.
We have
\begin{eqnarray*}
\int_{\bbbr^n} |u|\, d\mu
& \leq &
\sum_{k=-\infty}^\infty 2^k \mu(\{ 2^{k-1}<|u|\leq 2^k\}) \\
& \leq &
\sum_{k=-\infty}^\infty 2^k \mu(\{ |u|\geq 2^{k-1}\}) \\
& = &
\sum_{k=-\infty}^\infty 2^k \mu(\{ u_k\geq 2^{k-2}\})\\
& = &
\sum_{k=-\infty}^\infty 2^k \mu(\{ 2^{-(k-2)}u_k\geq 1\})\\
& \leq &
c(n)M \sum_{k=-\infty}^\infty 2^k \int_{\bbbr^n} |\nabla(2^{-(k-2)}u_k)|\, dx \\
& = & 
c(n)M \sum_{k=-\infty}^\infty 4\int_{\{2^{k-2}\leq |u|\leq 2^{k-1}\}} |\nabla u|\, dx\\
& = &
c(n)M\int_{\bbbr^n} |\nabla u|\, dx.
\end{eqnarray*}
The proof is complete.
\hfill $\Box$

\end{document}